\documentclass[11pt,twoside]{amsart}
\usepackage[latin1]{inputenc}
\usepackage[T1]{fontenc}
\usepackage{amsmath}
\usepackage{amsthm}
\usepackage{amssymb}
\usepackage[all]{xypic}

\newtheorem{thm}{Theorem}[section]

\newtheorem{prop}[thm]{Proposition}

\newtheorem{lem}[thm]{Lemma}
\newtheorem{cor}[thm]{Corollary}
\newtheorem{conj}[thm]{Conjecture}

\numberwithin{equation}{section}

\theoremstyle{definition}

\newtheorem{remark}[thm]{Remark}

\newcommand{\qqed}{\hspace*{\fill}$\Box$}
\newcommand{\supp}{\operatorname{Supp}}

\newcommand{\Db}{{\rm D}^{\rm b}}

\newcommand{\Aut}{{\rm Aut}}

\newcommand{\Pic}{{\rm Pic}}

\newcommand{\rk}{{\rm rk}}

\newcommand{\Spec}{{\rm Spec}}
\newcommand{\Spf}{{\rm Spf}}

\newcommand{\tr}{{\rm tr}}

\newcommand{\dual}{^{\vee}}

\newcommand{\cat}[1]{\begin{bf}#1\end{bf}}
\newcommand{\Ext}{{\rm Ext}}

\newcommand{\cal}{\mathcal}

\newcommand{\kc}{{\cal C}}

\newcommand{\ke}{{\cal E}}
\newcommand{\kf}{{\cal F}}

\newcommand{\kh}{{\cal H}}

\newcommand{\kn}{{\cal N}}
\newcommand{\km}{{\cal M}}
\newcommand{\ko}{{\cal O}}
\newcommand{\kp}{{\cal P}}

\newcommand{\kt}{{\cal T}}
\newcommand{\ku}{{\cal U}}
\newcommand{\kx}{{\cal X}}

\newcommand{\ZZ}{\mathbb{Z}}
\newcommand{\QQ}{\mathbb{Q}}
\newcommand{\RR}{\mathbb{R}}
\newcommand{\CC}{\mathbb{C}}
\newcommand{\EE}{\mathbb{E}}

\newcommand{\PP}{\mathbb{P}}

\newcommand{\Ch}{{\rm CH}}

\renewcommand{\to}{\xymatrix@1@=15pt{\ar[r]&}}
\renewcommand{\rightarrow}{\xymatrix@1@=15pt{\ar[r]&}}
\renewcommand{\mapsto}{\xymatrix@1@=15pt{\ar@{|->}[r]&}}
\renewcommand{\twoheadrightarrow}{\xymatrix@1@=15pt{\ar@{->>}[r]&}}
\renewcommand{\hookrightarrow}{\xymatrix@1@=15pt{\ar@{^(->}[r]&}}
\newcommand{\congpf}{\xymatrix@1@=15pt{\ar[r]^-\sim&}}
\renewcommand{\cong}{\simeq}
\newcommand{\Mod}[1]{{\ko_{#1}\text{-}\cat{Mod}}}

\begin{document}

\title[Chow groups of K3 surfaces and spherical objects]{Chow groups of K3 surfaces and spherical objects}
\author[D.\ Huybrechts]{Daniel Huybrechts}

\address{Mathematisches Institut,
Universit{\"a}t Bonn, Beringstr.\ 1, 53115 Bonn, Germany}
\email{huybrech@math.uni-bonn.de}

\begin{abstract} \noindent
We show that for a K3 surface $X$ the finitely generated subring
$R(X)\subset\Ch^*(X)$ introduced by Beauville and Voisin is
preserved under derived equivalences. This is proved by analyzing
Chern characters of spherical bundles (and complexes). As for a K3
surface $X$ defined over a number field all spherical bundles on
the complex K3 surface $X_\CC$ are defined over $\bar\QQ$, this is
compatible with the Bloch--Beilinson conjecture. Besides the work
of Beauville and Voisin \cite{BV}, Lazarfeld's result on
Brill--Noether theory for curves in K3 surfaces \cite{Laz} and the
deformation theory developed in \cite{HMS} are central for the
discussion. \vspace{-2mm}\end{abstract}

\maketitle
\let\thefootnote\relax\footnotetext{This work was supported by the SFB/TR 45 `Periods,
Moduli Spaces and Arithmetic of Algebraic Varieties' of the DFG
(German Research Foundation)}

\section{Introduction}

The Chow group $\Ch^i(X)$ of a smooth projective variety $X$ is
the group of all cycles of codimension $i$ modulo rational
equivalence (see \cite{Fulton}). For surfaces these are
$\Ch^0(X)=\ZZ[X]$, $\Ch^1(X)\cong\Pic(X)$ (via the first Chern
class), and the more mysterious $\Ch^2(X)$. The latter is roughly
the group of $0$-cycles $Z=\sum n_ix_i$ modulo linear equivalence
on curves containing $Z$. Since Mumford's article \cite{Mum} one
knows that, contrary to $\Ch^0(X)$ and $\Ch^1(X)$, the group
$\Ch^2(X)$ can be big. More precisely, the subgroup
$A(X):=\ker(\deg:\Ch^2(X)\to\ZZ)$ of all homologically trivial
$0$-cycles on a complex projective surface $X$ is
infinite-dimensional (and in particular infinitely generated)
whenever $p_g(X)=h^2(X,\ko_X)>0$, e.g.\ for K3 surfaces.

For K3 surfaces, Beauville and Voisin have studied more recently
the subgroup
$$R(X):=\Ch^0(X)\oplus\Ch^1(X)\oplus
c_X\ZZ\subset\Ch^2(X),$$ where $c_X\in\Ch^2(X)$ is the fundamental
class of a closed point $x\in X$ that is contained in a (possibly)
singular rational curve in $X$. As shown in \cite{BV}, the class
$c_X$ is independent of the point $x$. The main results of
\cite{BV} can be stated as follows:

\medskip

\noindent {\bf Theorem (Beauville--Voisin).} {\it {\rm i)}
$R(X)\subset \Ch^*(X)$ is a subring and {\rm ii)} $24c_X={\rm
c}_2(X)$.}

\medskip

The first condition is equivalently expressed by saying that for
any line bundle $L\in\Pic(X)$ the class ${\rm
c}_1(L)^2\in\Ch^2(X)$ is a multiple of $c_X$.

Let us rephrase i) and ii) in terms of Mukai vectors. For any
coherent sheaf (or complex of coherent sheaves) $E$ on $X$ one
defines
$$v^\Ch(E):={\rm ch}(E)\sqrt{{\rm td}(X)}\in\Ch^*(X).$$
In other words, $v^\Ch(X):=(\rk(E),{\rm c}_1(E),{\rm
ch}_2(E)+\rk(E)c_X)\in\Ch^0(X)\oplus\Ch^1(X)\oplus\Ch^2(X)$. Then
by applying i) to powers of $L$, conditions i) and ii) can be
reformulated as:

\begin{itemize}
  \item[iii)] {\it For any line bundle $L\in\Pic(X)$, one has
  $v^\Ch(L)\in R(X)$.}
\end{itemize}

\medskip

Line bundles on a K3 surface $X$ are the easiest examples of
\emph{spherical objects} on $X$, which by definition are bounded
complexes of coherent sheaves $E\in\Db(X)$ with
$\Ext^*_X(E,E)\cong H^*({\rm S}^2,\CC)$. Building upon \cite{BV},
we shall prove the following generalization of iii).

\medskip

\noindent{\bf Theorem 1.} {\it Let $X$ be a complex projective K3
surface of Picard number $\rho(X)\geq2$ and let $E\in\Db(X)$ be a
spherical object. Then $v^\Ch(E)\in R(X)$.}

\medskip

Spherical objects play a distinguished role in the study of the
bounded derived category $\Db(X)$ (and its homological mirror
given by a certain Fukaya category). They are essential for the
understanding of the rich structure of the group ${\rm
Aut}(\Db(X))$ of all exact $\CC$-linear autoequivalences and
Bridgeland's space of stability conditions ${\rm Stab}(X)$ (see
\cite{Br}).

In this context, Theorem 1 is used to deduce information about the
action of derived equi\-valences on the level of Chow groups. More
precisely, we have

\medskip

\noindent {\bf Theorem 2.} {\it Let
$\Phi_\ke:\Db(X)\congpf\Db(X')$ be an exact $\CC$-linear
equivalence between the bounded derived categories of two smooth
complex projective K3 surfaces of Picard number $\rho(X)\geq2$.
Then the induced action $\Phi_\ke^\Ch:\Ch^*(X)\congpf\Ch^*(X')$
preserves the Beauville--Voisin ring, i.e.\
$$\Phi^\Ch_\ke(R(X))=R(X').$$}

\medskip

The key step towards Theorem 1 is the following result which is
valid without any assumption on the Picard group of the surface
(see Theorem \ref{prop_consequencedefo}).

\medskip

\noindent {\bf Theorem 3.} {\it Let $X$ and $X'$ be complex
projective K3 surfaces and
$\Phi_\ke,\Phi_\kf:\Db(X)\congpf\Db(X')$ be two Fourier--Mukai
equivalences. If their induced actions on cohomology coincide,
i.e.\ $\Phi_\ke^H=\Phi_\kf^H:\widetilde H(X,\ZZ)\congpf \widetilde
H(X',\ZZ)$, then also $\Phi_\ke^{\rm CH}=\Phi_\kf^{\rm CH}:{\rm
CH}(X)\congpf{\rm CH}(X')$.}

\medskip

In particular, for $X= X'$ the result shows  that $\Phi_\ke^\Ch$
acts as the identity on $A(X)$ whenever it acts trivially on
cohomology (see Corollary \ref{cor:kernelsequal}). As we will
explain in Remark \ref{rem:firstBloch}, this is predicted by a
general conjecture of Bloch which asserts that the action of any
algebraic correspondence on the graded pieces of his conjectural
filtration is  determined by its action on cohomology (see
\cite[1.8]{Bl} or \cite[Conj.\ 23.22]{Voisin}).

\bigskip

If instead of  projective K3 surfaces over $\CC$ we consider  K3
surfaces $X$ defined over a number field $K$, then the situation
changes completely. In this case $\Ch^2(X)$ is no longer expected
to be infinitely generated. In fact, the Bloch--Beilinson
conjectures predict that for K3 surfaces over number fields the
degree map yields an isomorphism $\Ch^2(X)\otimes\QQ\cong\QQ$ (see
\cite{BlochCrelle,RSS}). This can be rephrased as the following

\bigskip

\noindent {\bf Conjecture (Bloch--Beilinson for K3 surfaces).}
{\it For any K3 surface $X$ defined over a number field $K$ base
change yields
$$\Ch^*(X)\otimes\QQ~\hookrightarrow~\Ch^*(X_L)\otimes\QQ\congpf
R(X_\CC)\otimes\QQ\cong\QQ^{\rho(X_\CC)+2},$$ where $L/K$ is a
certain finite field extension with a chosen embedding
$L\subset\CC$.}

\bigskip

As usual, $X_L$ denotes $X\times_{\Spec(K)}\Spec(L)$ and similarly
for $X_\CC$, which then is a complex projective K3 surface.

It is well known that the base change $\Ch^*(X)\to\Ch^*(X_L)$ for
any extension $L/K$ has torsion kernel.  The passage to the finite
extension $L/K$ is not essential and  only needed to ensure that
all geometric line bundles are defined. Thus, the central point of
the conjecture is that $\Ch^2(X)\otimes\QQ\cong\QQ c_X$. For a
proof it would clearly suffice to prove that any rational point
$x\in X(\bar\QQ)$ satisfies $[x]=c_X$, but there is no obvious
geometric reason for this.

 The skyscraper sheaves $k(x)$ of rational points $x\in X(\bar\QQ)$ define
semi-rigid objects in $\Db(X_{\bar\QQ})$ (see Section
\ref{sect:K3overnumber} for the definition), for
$\Ext^1(k(x),k(x))$ is two-dimensional. In this sense, they are
reasonably close to our spherical objects $E\in\Db(X)$ which have
vanishing $\Ext^1(E,E)$. The following is thus in accordance with
the Bloch--Beilinson conjecture for K3 surfaces.

\medskip

\noindent {\bf Theorem 3.} {\it Let $X$ be a K3 surface over a
number field $K$. Then

$\bullet$ Any spherical object $E\in\Db(X_\CC)$ is defined over
some finite extension $L/K$.

$\bullet$ Any spherical object $F\in\Db(X_{\bar\QQ})$ satisfies
$v^\Ch(F_\CC)\in R(X_\CC)$ if in addition
$\rho(X_{\bar\QQ})\geq2$.}

\bigskip

Together these two assertions show that on any K3 surface $X$
defined over $\bar\QQ$ with Picard number $\rho(X)\geq2$ there
exists a large number of non-trivial classes in $\Ch^2(X)$ for
which the Bloch--Beilinson conjecture can be verified. As with
other approaches to the conjecture, the difficult part, that would
show that this suffices to deduce the result for all classes,
remains open.

\medskip

Despite the algebraic nature of all assertions, non-algebraic K3
surfaces play a crucial but hidden role in this paper. A technique
that has been developed together with Macr\`i and Stellari in
\cite{HMS} allows one to deform any derived equivalence that acts
as the identity on cohomology to an essentially trivial derived
equivalence on a generic and non-algebraic deformation of $X$.
This is explained in Section \ref{sect:FMonChow}. The results can
also be used to show that the natural representations of ${\rm
Aut}(\Db(X))$ on $\Ch^*(X)$ and on the Mukai lattice $\widetilde
H(X,\ZZ)$ encode the same information (see Corollary
\ref{cor:kernelsequal}).

In Section \ref{sect:Mukaispherobjcts} we reduce Theorem 1 to the
case that $E$ is a spherical vector bundle. This section also
explains how Theorem 2 is deduced from Theorem 1. The case of
spherical vector bundles is dealt with in detail in Section
\ref{sect:MukaisphVB}. The main ingredient here is Lazarsfeld's
result that the generic curve in an indecomposable linear system
on a K3 surface is Brill--Noether general \cite{Laz}. The final
Section \ref{sect:K3overnumber} discusses the relation to the
Bloch--Beilinson conjecture for K3 surfaces over number fields.

We certainly expect all results to hold true without any
assumption on the Picard number of the K3 surface $X$ and, in
fact, Theorem 1 can be proved also for $\rho(X)=1$ under
additional numerical conditions on the spherical object $E$.

\bigskip

\noindent {\bf Notation.} By $\Db(X)$ we denote the bounded
derived category of the abelian category ${\rm Coh}(X)$ of
coherent sheaves on $X$. It will be considered as a $K$-linear
triangulated category, when $X$ is defined over $K$ (which mostly
is $\CC$, $\bar\QQ$ or a number field). The Mukai lattice
$\widetilde H(X,\ZZ)$ of a complex K3 surface is by definition the
full singular cohomology $H^*(X,\ZZ)$ endowed with its natural
weight two Hodge structure and the Mukai pairing (see e.g.\
\cite{FM}). All intersection products will be taken with respect
to the Mukai pairing which differs from the usual intersection
pairing by a sign in degree four. We will associate to any
$E\in\Db(X)$ its Mukai vector $v(E):={\rm ch}(E)\sqrt{{\rm
td}(X)}\in\widetilde H(X,\ZZ)$ and its natural lift to $\Ch^*(X)$,
which is denoted $v^\Ch(E)$. We do not make this distinction for
the characteristic classes.

\bigskip

\noindent {\bf Acknowledgements.} I wish to thank C.\ Voisin for
instructive comments on a first version of the paper and R.\
Kloosterman for his help with an argument in Section
\ref{sect:K3overnumber}. Thanks also to E.\ Mistretta, M.\
Penegini, and D.\ Ploog, who have checked Theorem 2 in explicit
examples at an early stage of this work.

\section{Fourier--Mukai action on the Chow
group}\label{sect:FMonChow}

Let us start by briefly recalling the following examples of
Fourier--Mukai equivalences for K3 surfaces.

\medskip

{\bf i)} For a line bundle $L\in\Pic(X)$, the tensor product
$L\otimes(~~)$ defines a Fourier--Mukai equivalence
$$\Phi_{\iota_*L}:\Db(X)\congpf\Db(X),$$
where $\iota_*L$ is the direct image of $L$ under the diagonal
embedding $X\congpf\Delta\subset X\times X$.

\smallskip

{\bf ii)} If $X'$ is a smooth projective two-dimensional fine
moduli space of $\mu$-stable vector bundles on $X$, then the
universal bundle $\EE$ on $X\times X'$ induces an equivalence (see
\cite{Muk,FM})
$$\Phi_\EE:\Db(X)\congpf\Db(X').$$
Note that any Fourier--Mukai partner of $X$ is isomorphic to such
a moduli space, but of course other Fourier--Mukai equivalences
between $X$ and $X'$ do exist and are given by kernels more
complicated than $\EE$.

\smallskip

{\bf iii)} If $E\in\Db(X)$ is a spherical object, i.e.\
$\Ext^*_X(E,E)\cong H^*({\rm S}^2,\CC)$, then the spherical twist
$$T_E:\Db(X)\congpf\Db(X),$$
studied in detail in \cite{ST}, is a Fourier--Mukai equivalence
with kernel $\kp_E:={\rm Cone}(\tr:E\dual\boxtimes
E\to\ko_\Delta)$ (see also \cite[Ch.\ 8]{FM}). Spherical objects,
although rigid, exist in abundance on any projective K3 surface.
E.g.\ any line bundle, even the trivial one, gives rise to a
non-trivial spherical twist. Moreover, Kuleshov shows in
\cite{Kul} that any $(1,1)$-class  $v\in\widetilde H(X,\ZZ)$ of
square $-2$ is the Mukai vector of a spherical object, which can
be chosen to be a vector bundle if the rank of $v$ is positive.

\medskip

Any Fourier--Mukai equivalence $\Phi_\ke:\Db(X)\congpf\Db(X')$
induces a group isomorphism
$$\Phi^\Ch_\ke:\Ch^*(X)\congpf\Ch^*(X')$$
and a Hodge isometry $$\Phi^H_\ke:\widetilde H(X,\ZZ)\congpf
\widetilde H(X',\ZZ).$$ Both are defined as correspondences
associated to $v^\Ch(\ke)\in \Ch^*(X\times X')\otimes\QQ$
respectively $v(\ke)\in H^*(X\times X',\QQ)$.

Note that in general one would expect $\Phi^\Ch$ and $\Phi^H$ to
be defined only with rational coefficients, but as Mukai observed
the situation is special for K3 surfaces (see the original
argument in \cite{Muk} or \cite{FM,HL}).

\begin{remark}
The action $\Phi^\Ch_\ke$ is difficult to grasp for example ii),
but easy to describe in the examples i) and iii).

Indeed, in i) the actions $\Phi^\Ch_{\iota_*L}$ and
$\Phi^H_{\iota_*L}$ are both given by multiplication with ${\rm
ch}(L)=\exp({\rm c}_1(L))$, where the Chern character is viewed in
$\Ch^*(X)$ resp.\ $H^*(X,\ZZ)$. Thus Theorem 2 is a trivial
consequence of the results in \cite{BV} in this case.

For  the spherical twists in iii), $T^\Ch_E$ and $T^H_E$ are
reflections in $v^\Ch(E)^\perp$ resp.\ $v(E)^\perp$, where the
orthogonal complement is taken with respect to the Mukai pairing.
In particular, their squares $(T^2_E)^\Ch$ and $(T^2_E)^H$ act
trivially, i.e.\ as the identity, on both groups $\Ch^*(X)$ resp.\
$\widetilde H(X,\ZZ)$.
\end{remark}

\begin{remark}\label{rem:actioninexamples}
Observe that for an arbitrary spherical object the associated
spherical twist $T_E$ preserves the Beauville--Voisin subring
$R(X)$ if and only if $v^\Ch(E)\in R(X)$. In this case it acts as
the identity on the space of cohomologically trivial cycles
$A(X)$.
\end{remark}

According to Mumford, $\Ch^2(X)$ is big and in fact of infinite
dimension for any complex projective surface with $p_g(X)>0$ and
therefore in particular for K3 surfaces. See \cite[Ch.\
22]{Voisin} for the notion of dimension of $\Ch^2(X)$.

Thus, a priori for an arbitrary Fourier--Mukai equivalence
$\Phi_\ke:\Db(X)\congpf\Db(X')$ between two K3 surfaces the
induced map $\Phi^\Ch_\ke:\Ch^*(X)\congpf\Ch^*(X')$ between the
infinite dimensional  Chow groups might capture more information
than $\Phi_\ke^H:\widetilde H(X,\ZZ)\congpf \widetilde H(X',\ZZ)$.
That this is (unfortunately?) not the case is the main result of
this section

\begin{thm}\label{prop_consequencedefo}
Let $X$ and $X'$ be smooth complex projective K3 surfaces and let
$$\Phi_\ke,\Phi_\kf:\Db(X)\congpf\Db(X')$$ be two Fourier--Mukai
equivalences with $\Phi_\ke^H=\Phi_\kf^H$. Then also
$$\Phi_\ke^\Ch=\Phi^\Ch_\kf:\Ch^*(X)\congpf \Ch^*(X').$$
\end{thm}

\begin{remark} In general the direct sum decomposition  $\Ch^*(X)=\Ch^0(X)\oplus\Ch^1(X)\oplus\Ch^2(X)$
is not respected by Fourier--Mukai transforms. However, the
homologically trivial part is. More precisely, if
$\Phi_\ke:\Db(X)\congpf\Db(X')$ is any Fourier--Mukai equivalence,
then $\Phi_\ke^\Ch(A(X))=A(X')$. Moreover, if $\Phi^H_\ke$
respects the cohomological degree, e.g.\ for cohomologically
trivial autoequivalences, then
$\Phi_\ke^\Ch(\Ch^0(X)\oplus\Ch^2(X))=\Ch^0(X')\oplus \Ch^2(X')$.
\end{remark}

The essential step in the proof of Theorem
\ref{prop_consequencedefo} consists of the following slightly
weaker result.

\begin{prop}\label{prop:weak}
Let $\Phi_\ke,\Phi_\kf:\Db(X)\congpf\Db(X')$ be as in Theorem
\ref{prop_consequencedefo}. Then $\Phi^\Ch_\ke=\Phi_\kf^\Ch$ on
$\Ch^0(X)\oplus\Ch^2(X)$.
\end{prop}

\begin{proof}
By studying the composition
$\Phi_\kf^{-1}\circ\Phi_\ke:\Db(X)\congpf\Db(X)$, one easily
reduces to the case of autoequivalences acting as the identity on
cohomology. So let $\Phi_{\ke_0}:\Db(X)\congpf\Db(X)$ with
$\Phi^H_{\ke_0}={\rm id}$. We claim that then also
$\Phi^\Ch_{\ke_0}={\rm id}$ on $\Ch^0(X)\oplus \Ch^2(X)$. In
particular we have to show that $\Phi^\Ch_{\ke_0}={\rm id}$ on the
space of homologically trivial cycles $A(X)$.

Clearly, changing $\Phi_{\ke_0}$ by even powers $T^{2k}$ of the
shift functor or even powers $T_{\ko_X}^{2 k}$ of the spherical
twist associated to the trivial line bundle does not affect the
assertion (see Remark \ref{rem:actioninexamples} and use
$v^\Ch(\ko_X)=(1,0,c_X)\in R(X)$). So, in the course of the proof
we will freely modify $\Phi_{\ke_0}$ by autoequivalences of this
type.

In \cite{HMS} we were mainly interested in the case
$\Phi_{\ke_0}^H=(-{\rm id}_{H^2})\oplus {\rm id}_{H^0\oplus H^4}$,
but as mentioned there already the case $\Phi_{\ke_0}^H={\rm id}$
is similar and actually easier. So the results of \cite{HMS}  show
that for any autoequivalence $\Phi_{\ke_0}$
with $\Phi^H_{\ke_0}={\rm id}$ one finds:\\

i) Two smooth formal deformations
$\kx\to\Spf(R)\xymatrix@1@=15pt{&\ar[l]}\kx'$ with $R=\CC[[t]]$
and $\kx_0=X=\kx'_0$. Here $\kx$ is the formal neighbourhood of
$X$ inside its twistor space with respect to a very general
K\"ahler class in $\Pic(X)\otimes\RR$. Note that in this way $X$
is deformed towards a non-projective K3 surface.

ii) A complex $\ke\in\Db(\kx\times_R\kx'):=\Db_{\rm
coh}(\Mod{\kx\times_R\kx'})$ deforming $\ke_0$, i.e.\
$L\iota^*\ke\cong\ke_0$, where
$\iota:X\times X ~\hookrightarrow \kx\times_R^{\phantom{R}}\kx'$ is the obvious closed embedding.\\

By \cite[Prop.\ 2.18, 2.19]{HMS} we may assume, after possibly
composing with powers of $T_\ko^2$ and $T^2$, that the restriction
$\ke_K\in\Db((\kx\times_R\kx')_K)$ of $\ke$ to the general fibre
is a sheaf. Hence \cite[Cor.\ 4.5]{HMS} applies and shows that
there exists an $R$-flat sheaf(!) $\widetilde\ke$ on
$\kx\times_R\kx'$ with the same restriction to the general fibre
as $\ke$, i.e.\ $\widetilde\ke_K\cong\ke_K$ in
$\Db((\kx\times_R\kx')_K)$, where $K=\CC((t))$. For the notation
we refer to \cite{HMS}. Using the compatibilities between
$\Phi_{\ke_0}^H$ and the induced action on Hochschild
(co)homology, one can in addition assume that the first order
deformations
$\kx_1\to\Spec(\CC[t]/t^2)\xymatrix@1@=15pt{&\ar[l]}\kx'_1$ of
$X=\kx_0=\kx'_0$ coincide.

The specialization morphism $K(\Db((\kx\times_R\kx')_K))\to
K(\Db(X\times X))$ is well defined, see \cite[Remark 2.7]{HMS} or
the analogous statement for Chow groups in Remark \ref{rem:spec}.
Hence the coherent sheaf(!) $\widetilde\ke_0$ and the original
complex $\ke_0$ have the same Mukai vectors $v^\Ch\in\Ch^*(X\times
X)\otimes\QQ$ and therefore
$\Phi_{\widetilde\ke_0}^\Ch=\Phi_{\ke_0}^\Ch$ and
$\Phi_{\widetilde\ke_0}^H=\Phi_{\ke_0}^H$.

Note that the Fourier--Mukai transform $\Phi_G$ associated to the
sheaf $G:=\widetilde\ke_0$ is not necessarily an equivalence,
which would simplify the following arguments. But in any case,
there is a dense open subset $U\subset X$ over which $G$ is flat
(see e.g.\ \cite[Thm.\ 2.15, Lemma 2.1.6]{HL}). Hence, for any
closed point $x\in U$ the image $\Phi_G(k(x))$ is simply the sheaf
$G|_{\{x\}\times X}$. On the other hand,
$v(\Phi_G(k(x)))=v(\Phi_{\ke_0}(k(x)))=(0,0,1)$ and hence
$G|_{\{x\}\times X}$ must be of the form $k(y)$ for some point
$y\in X$. This gives rise to a morphism $U\to X$, which by
interchanging the two factors turns out to define a birational map
$X\dashrightarrow X$. As any birational map between K3 surfaces,
the latter can then be completed to an isomorphism $f:X\congpf X$.
Moreover, if $Z:=\supp(G)\subset X\times X$, then $\Gamma_f\subset
Z$ is one irreducible component and the other components do not
dominate $X$. The latter implies that
$[\Gamma_f]_*|_{H^{2,0}}=\Phi^H_G|_{H^{2,0}}={\rm id}_{H^{2,0}}$,
i.e.\ $f$ is a symplectomorphism.

Obviously $f_*(c_X)=c_X$ and a general conjecture of Bloch (see
Remark \ref{rem:firstBloch}) predicts that for a symplectomorphism
the induced automorphism $f_*:\Ch^*(X)\cong\Ch^*(X)$ is the
identity on $A(X)$. Since for generic $x\in X$ we have
$\Phi_G^\Ch([x])=[f(x)]$, this would be enough to conclude that
$\Phi_G^\Ch$ acts as the identity on $\Ch^2(X)$.

Without using Bloch's conjecture, the argument is more involved
and goes as follows. Since $G$ is the restriction of a sheaf on
$\kx_1\times_{R_1}\kx_1$, the structure sheaf of the graph
$\ko_{\Gamma_f}$ deforms sideways to first order, i.e.\ there
exists an $R_1$-flat coherent sheaf on $\kx_1\times_{R_1}\kx'_1$
restricting to $\ko_{\Gamma_f}$ over the closed point. (Do it
first for the graph of $f|_U$ and then pass to the closure.) In
other words, the automorphism $f$ deforms sideways to first order
(actually to any order, but we do not need this) in
$\kx\times_{R}\kx'$. But clearly $f$ deforms sideways to first
order if and only if $f^*(w)=w$, where $w\in H^1(X,\kt_X)$
corresponds to the first order deformation
$\kx_1\to\Spec(\CC[t]/t^2)$.

By construction, the class $w$ maps to the chosen  K\"ahler class
in $H^{1,1}(X)$ under the isomorphism $H^{1,1}(X)\cong
H^1(X,\Omega_X)\cong H^1(X,\kt_X)$ and, since the K\"ahler class
was chosen generically, this implies $f^*={\rm id}$ on $\Pic(X)$.
Since the transcendental lattice $T(X)$ is an irreducible Hodge
structure (of weight two), the assumption $f^*={\rm id}$ on
$H^{2,0}(X)$ implies by Schur's lemma that $f^*={\rm id}$ on
$T(X)$. Together with $f^*={\rm id}$ on ${\rm Pic}(X)$ this proves
$f^*={\rm id}$ on the full cohomology $H^*(X,\ZZ)$. By the Global
Torelli theorem, the latter is equivalent to $f={\rm id}$.
Eventually this shows that $\Phi_G^\Ch([x])=f_*[x]=[x]$ for
generic and hence all $x\in X$. Thus $\Phi^\Ch_G={\rm id}$ on
$\Ch^2(X)$.

To conclude we observe that $\Phi_{\ke_0}(\ko_X)$ deforms sideways
to a spherical object in $\Db(\kx_K)$, for $\ko_X$ and $\ke_0$ do.
On the other hand, up to shift $\ko_{\kx'_K}$ is the only
spherical object in $\Db(\kx'_K)$ (cf.\ \cite[Prop.\ 2.14]{HMS}).
Hence, up to shift $\ko_X$ is the only spherical object on $X$
that deforms sideways in the family $\kx'$ to a spherical object
in $\Db(\kx_K')$. Hence $\Phi_{\ke_0}(\ko_X)\cong\ko_X$ (up to
shift), which in particular shows that
$\Phi^\Ch_{\ke_0}(1,0,c_X)=(1,0,c_X)$.

Thus we have shown that $\Phi_{\ke_0}^\Ch$ acts as identity on
$\Ch^0(X)\oplus\Ch^2(X)$. Moreover, it acts as $\left(\begin{array}{cc}{\rm id}&0\\
\ast&{\rm id}\end{array}\right)$ on $\Ch^1(X)\oplus\Ch^2(X)$ which
will later be shown to be diagonal.
\end{proof}

\medskip

The proposition can also be used to derive information about the
Mukai vectors in $\Ch^*(X)$ of spherical objects having the same
Mukai vector in cohomology. This is the following

\begin{cor}\label{cor:sphsame}
If $E,E'\in\Db(X)$ are two spherical objects with $v(E)=v(E')\in
\widetilde H(X,\ZZ)$, then $$v^\Ch(E)=v^\Ch(E')\in \Ch^*(X).$$
\end{cor}

\begin{proof}
Write $v(E)=(r,\ell,s)=v(E')$. Let us first reduce to the case
that $r\ne0$. Suppose $r=0$, then $\ell\ne0$. Then let
$\Phi=T_{\ko_X}\circ (L\otimes(~~))$ and use $\Phi\circ T_E\cong
T_{\Phi(E)}\circ \Phi$ (see e.g.\ \cite[Lemma 8.21]{FM}), which
holds for any Fourier--Mukai equivalence $\Phi$ and any spherical
object $E$. Thus $T^H_{\Phi(E)}=T^H_{\Phi(E')}$ and the assertion
$v^\Ch(E)=v^\Ch(E')$ is clearly equivalent to
$v^\Ch(\Phi(E))=v^\Ch(\Phi(E'))$. If  $L$ is chosen such that
$({\rm c}_1(L).\ell)\ll 0$, then the spherical object $\Phi(E)$
has positive rank.

Now apply Proposition \ref{prop:weak} to the class $c_X$ to deduce
 $T^\Ch_E(c_X)=T^\Ch_{E'}(c_X)$. Both sides can be explicitly
computed, which yields $c_X-rv^\Ch(E)=c_X-rv^\Ch(E')$ and hence
$v^\Ch(E)=v^\Ch(E')$.
\end{proof}

\noindent {\it Proof of Theorem \ref{prop_consequencedefo}.}
Suppose again that $\Phi_\ke$ is an autoequivalence of $\Db(X)$
with $\Phi_\ke^H={\rm id}$. By Proposition \ref{prop:weak} we know
already that $\Phi_\ke^\Ch$ is the identity on
$\Ch^0(X)\oplus\Ch^2(X)$. Thus it remains to show that
$\Phi^\Ch_\ke({\rm c}_1(L))$ of an arbitrary line bundle $L$ has
no component in $A(X)$. The image $L':=\Phi_\ke(L)$ of a line
bundle $L$ is a spherical object and since $\Phi^H_\ke={\rm id}$,
one has $v(L)=v(L')$. By Corollary \ref{cor:sphsame} this implies
$v^\Ch(L)=v^\Ch(L')$ and hence $\Phi^\Ch_\ke({\rm c}_1(L))\in
\Ch^1(X)$. \qqed

\bigskip

As done already in the proof above, Theorem
\ref{prop_consequencedefo} can be reformulated in terms of
autoequivalences. Since the kernel of the cohomology
representation of ${\rm Aut}(\Db(X))$ is essentially the only
remaining  mystery in this context, we state this explicitly as

\begin{cor}\label{cor:kernelsequal}
Let $X$ be a smooth complex projective K3 surface and denote by
$$\rho^\Ch:{\rm Aut}(\Db(X))\to{\rm Aut}(\Ch^*(X))\phantom{~~~~~}{and}~~~~~ \rho^H:{\rm
Aut}(\Db(X))\to{\rm Aut}(\widetilde H(X,\ZZ))$$ the natural
representation $\Phi\mapsto\Phi^\Ch$ resp.\ $\Phi\mapsto\Phi^H$.
Then $\ker(\rho^\Ch)=\ker(\rho^H)$.
\end{cor}

\begin{proof}
The inclusion $\ker (\rho^\Ch)\subset\ker(\rho^H)$ is obvious and
the other one follows from the proposition.
\end{proof}

\begin{remark} In \cite{Br}
Bridgeland suggests the following explicit description of this
kernel. He conjectures  $\ker(\rho^H)=\pi_1(\kp_0^+(X))$, for a
certain period domain $\kp^+_0(X)$ defined in terms of the
algebraic part of $\widetilde H(X,\ZZ)$. In particular, the
conjecture says that $\ker(\rho^H)$ is spanned by the square
$T^2=[2]$ of the shift functor and the squares $T^2_E$ of all
spherical twists $T_E$. (In fact, spherical twists associated to
spherical sheaves should suffice.)

As explained above, the conjectural generators $T^2$ and $T_E^2$
of $\ker(\rho^H)$ act trivially on $\Ch^*(X)$. In this sense, the
corollary provides non-trivial evidence for Bridgeland's
conjecture.
\end{remark}

\begin{remark}\label{rem:firstBloch}
The corollary appears interesting also in the light of another
open conjecture due to Bloch (see \cite{Bl}), which for the case
of surfaces  reads as follows: Consider a surface $X$ and a cycle
$\Gamma\in\Ch^2(X\times X)$ with its induced natural endomorphisms
$[\Gamma]^{2,0}_*$ of $H^0(X,\Omega_X^2)$ and $[\Gamma]_*$ of
$\Ch^2(X)$. In general the latter does not respect the natural
filtration $\ker({\rm alb}_X)\subset A(X)\subset\Ch^2(X)$, but
induces an endomorphism ${\rm gr}[\Gamma]_*$ of the graded object
$\ker({\rm alb}_X)\oplus{\rm Alb}(X)\oplus\ZZ$. Then Bloch
conjectures that $[\Gamma]_*^{2,0}=0$ if and only if ${\rm
gr}[\Gamma]_*$ is trivial on $\ker({\rm alb}_X)$ (see also
\cite[Ch.\ 11]{Voisin}). Note that for a K3 surface $X$ the graded
object is just $A(X)\oplus\ZZ$.

If $\Phi_\ke:\Db(X)\congpf\Db(X)$ is a Fourier--Mukai
autoequivalence of a K3 surface $X$ such that $\Phi^H_\ke={\rm
id}$, then $\Gamma:=v^\Ch(\ke)-[\Delta]$ acts trivially on
cohomology and, in particular, on $H^{0}(X,\Omega_X^2)$. Bloch's
conjecture would thus say that ${\rm gr}[\Gamma]_*$ is trivial on
$A(X)=\ker({\rm alb}_X)$. And indeed, by Corollary
\ref{cor:kernelsequal}  this holds true, as we in fact have
$\Gamma=0$.

Note that Bloch's conjecture would actually say that
$\Phi_\ke^H={\rm id}$ on $H^0(X,\Omega_X^2)$ is sufficient  to
conclude $\Phi_\ke^\Ch={\rm id}$ on $A(X)$, but our techniques
fail to prove this. In fact, it seems even unknown whether any
symplectomorphism $f\in \Aut(X)$ induces the identity on $A(X)$,
i.e.\ whether  for a symplectomorphism $f$ any point $x\in X$ is
rationally equivalent to its image $f(x)$ (cf.\ comments in the
proof of Proposition \ref{prop:weak}).
\end{remark}

\section{Mukai vectors of spherical
objects}\label{sect:Mukaispherobjcts}

The goal of this section is to prove that for any spherical object
$E\in\Db(X)$ on a smooth complex projective K3 surface $X$ with
Picard number $\rho(X)\geq2$ the Mukai vector
$v^\Ch(E)\in\Ch^*(X)$ is contained in the Beauville--Voisin
subring $R(X)\subset\Ch^*(X)$. The main results (Corollaries
\ref{cor:sphobR} and \ref{cor:2}) should also hold for K3 surfaces
with $\rho(X)=1$, but we can only prove it under additional
conditions on the numerical invariants of $E$.

\begin{remark}\label{rem:equivcondiandii}
As noted earlier, the original result in \cite{BV} can be seen as
the special case that the spherical object $E$ is a line bundle
$L\in\Pic(X)$. Indeed, if $v^\Ch(L^k)=(1,k{\rm c}_1(L),k^2{\rm
c}_1(L)^2/2+{\rm c}_2(X)/24)\in R(X)$, then necessarily ${\rm
c}_1(L)^2\in R(X)$ and ${\rm c}_2(X)\in R(X)$. However, it should
be emphasized that our methods do not provide an alternate proof
of the results in \cite{BV}.
\end{remark}

In this section we shall explain how to reduce the proof of
Theorem 1 to a generalization of the Beauville--Voisin result from
line bundles to higher rank spherical vector bundles. The proof of
the following crucial result is postponed to the next section.

\begin{prop}\label{prop:sphvb}
Let $X$ be a smooth complex projective K3 surface and let $E$ be a
spherical vector bundle on $X$. Suppose that one of the following
conditions hold \begin{itemize}\item[i)] $\rho(X)\geq2$ or
\item[ii)] $\Pic(X)=\ZZ H$ and $v(E)=(r,kH,s)$ with $k\equiv\pm1
(r)$.
\end{itemize}
Then $v^\Ch(E)\in R(X)$.
\end{prop}

This proposition is expected to hold without any restriction on
$X$ or the spherical vector bundle $E$. However, it does not
generalize to $\mu$-stable vector bundles, i.e.\ the Mukai vector
of a general non-rigid $\mu$-stable vector bundle $E$ is certainly
not contained in $R(X)$.

The following consequence of Proposition \ref{prop:sphvb}  proves
Theorem 1.


\begin{cor}\label{cor:sphobR}
Let $E\in\Db(X)$ be a spherical object on a smooth projective K3
surface $X$. Suppose that either
\begin{itemize}
\item[i)] $\rho(X)\geq2$ or \item[ii)] $\Pic(X)=\ZZ H$ with
$v(E)=(r,kH,s)$ with $k\equiv\pm1(r)$.
\end{itemize}
Then $v^\Ch(E)\in R(X)$.
\end{cor}

\begin{proof}
Let $E\in\Db(X)$ be spherical. Write $v(E)=(r,\ell,s)$. Clearly,
$v^\Ch(E)\in R(X)$ is equivalent to $v^\Ch(E[1])=-v^\Ch(E)\in
R(X)$. Hence we may assume $r\geq0$.

If $r>0$, then as proved in \cite{Kul} there exists a spherical
locally free sheaf $E'$ on $X$ with $v(E')=v(E)$. (In fact any
torsion free spherical sheaf is automatically locally free as was
already observed by Mukai in \cite{Muk}.) Now apply Proposition
\ref{prop:sphvb} to $E'$ which yields $v^\Ch(E')\in R(X)$. But by
Corollary \ref{cor:sphsame}, we know that $v^\Ch(E)=v^\Ch(E')$ for
any two numerically equivalent spherical objects $E,E'\in\Db(X)$.

The case $r=0$ is straightforward. First, by applying $T_{\ko_X}$
we reduce to the case that also $s=0$ and hence $v(E)=(0,\ell,0)$
with $\ell$ a $(-2)$-class, which we may assume to be effective.
Thus $T_E^H=s_\ell$, the reflection  in $\ell^\perp$. The Weyl
group $W_X$ generated by reflections $s_\delta$ for all effective
$(-2)$-classes $\delta$ is known to be generated by reflections
associated to nodal classes, i.e.\ when $\delta$ is represented by
a smooth rational curve. Thus $s_\ell$ can be written as a
composition of
 finitely many reflections $s_{[C_i]}$, where
the curves $C_i\subset X$ are smooth and rational. Then use that
$v^\Ch(E)\in R(X)$ is equivalent to $T^\Ch_E(R(X))=R(X)$, because
$T^\Ch_E$ is the reflection in the hyperplane  $v^\Ch(E)^\perp$,
and that $T^H_E=s_\ell=\prod s_{[C_i]}$ implies $T^\Ch_E=\prod
T_{\ko_{C_i}(-1)}^\Ch$ (Theorem \ref{prop_consequencedefo}).
Clearly, $v^\Ch(\ko_{C_i}(-1))\in R(X)$ and hence
$T^\Ch_{\ko_{C_i}(-1)}(R(X))=R(X)$.
\end{proof}

Let us now show that Theorem 1 implies Theorem 2, which we state
again as

\begin{cor}\label{cor:2}
Suppose $X$ and $X'$ are smooth complex projective K3 surfaces
with Picard number $\rho(X)\geq2$. If
$\Phi_\ke:\Db(X)\congpf\Db(X')$ is a Fourier--Mukai equivalence,
then the induced map $\Phi_\ke^\Ch:\Ch^*(X)\congpf\Ch^*(X')$
respects the Beauville--Voisin subring, i.e.\
$$\Phi_\ke^\Ch(R(X))=R(X').$$
\end{cor}

\begin{proof}
The Mukai vectors $v^\Ch(L)\in\Ch^*(X)$ of all line bundles
$L\in\Pic(X)$ span $R(X)$. The images $E:=\Phi_\ke(L)\in\Db(X')$
are not necessarily (shifted) line bundle again, but they are
spherical objects in $\Db(X')$. Since $\rho(X)\geq2$ and hence
$\rho(X')\geq2$, Corollary \ref{cor:sphobR}, i) applies. Thus
$\Phi_\ke^\Ch(v^\Ch(L))=v^\Ch(E)\in R(X')$.
\end{proof}

\begin{remark}
The most interesting special case is the one when $X'$ is a fine
moduli space of $\mu$-stable vector bundles and $\ke$ is the
universal bundle $\EE$. In examples where both, moduli space $X'$
and universal bundle $\EE$, are constructed explicitly, one
sometimes can prove Corollary \ref{cor:2} directly (see e.g.\
\cite{HL} for explicit examples). If there was an argument proving
the result for arbitrary universal bundles without first proving
Corollary \ref{cor:sphobR}, then the techniques of Section
\ref{sect:FMonChow} would prove Corollary \ref{cor:2} more
directly (and also in the case $\rho(X)=1$).

To be more precise, let $E$ be any spherical object with
$v(E)=(r,\ell,s)$ and $r>0$. Then
$T^H_E(0,0,1)=-(r^2,r\ell,rs-1)$. Thus, if $T_E$ is composed with
$\Phi_{\EE[1]}$ where $\EE\in\Db(X\times X')$ is the universal
family of stable vector bundles with Mukai vector
$(r^2,r\ell,rs-1)$, then $(\Phi_{\EE[1]}\circ
T_E)^H(0,0,1)=(0,0,1)$. By composing with a certain equivalence
$\Psi$ that is a combination of spherical twists $T_{\ko_C}$ (with
$\PP^1\cong C\subset X$) and tensor products with line bundles,
the Hodge isometry $(\Psi\circ\Phi_{\EE[1]}\circ T_E)^H$ becomes
graded. Moreover, it will respect the K\"ahler cone up to sign
(see e.g.\ \cite[Ch.\ 9]{FM} for details). Then by the Global
Torelli theorem $(\Psi\circ\Phi_{\EE[1]}\circ T_E)^H=\pm f_*$ for
some isomorphism $f$. Now use that $f^\Ch_*$ and $\Psi^\Ch$
preserve the Beauville--Voisin ring. Hence $T_E^\Ch(R(X))=R(X)$
 if and only if $\Phi^\Ch_{\EE}(R(X))=R(X')$, where  $\EE$
is a universal family of stable bundles of rank $r^2$.
\end{remark}



\section{Spherical vector bundles: Proof of Proposition
\ref{prop:sphvb}}\label{sect:MukaisphVB}

Let $C$ be a smooth irreducible complex projective curve of genus
$g$. Recall that the Brill--Noether locus
$W^{r_0}_d(C)\subset\Pic^d(C)$ is the determinantal subvariety of
all line bundles $A$ of degree $d$ with $h^0(C,A)\geq r_0+1$. The
Brill--Noether number for these numerical invariants is by
definition
$$\rho(r_0,d,g):=g-(r_0+1)(g-d+r_0).$$

Classically (see \cite{ACGH}) one knows that $W^{r_0}_d(C)$ is
non-empty whenever $\rho(r_0,d,g)\geq0$. (Due to a result of
Fulton and Lazarsfeld, it is also connected when
$\rho(r_0,d,g)>0$, but this will not be used.) Moreover, for a
generic curve $C$ the Brill--Noether number $\rho(r_0,d,g)$ is in
fact the dimension of $W^{r_0}_d(C)$ when $\rho(r_0,d,g)\geq0$ and
$W^{r_0}_d(C)=\varnothing$ otherwise.

Central for  our discussion  is a result of Lazarsfeld \cite{Laz}
that shows that a generic smooth curve $C$ in an indecomposable
linear system on a K3 surface is Brill--Noether general, i.e.\ the
$W^{r_0}_d(C)$ have the expected dimension. Let us make precise
which parts of \cite{Laz} are really used.

\medskip

Suppose $A\in W^{r_0}_d(C)$ satisfies
\begin{equation}\label{eqn:Lazcond1}
{\rm i)}~~h^0(C,A)=r_0+1 ~~{\rm ~and~ii)~the~line~bundles}~~
A~~{\rm and}~~ A^*\otimes\omega_C~~{\rm are~ globally~ generated}.
\end{equation}

If $C$ is embedded into a K3 surface $X$, one associates to $A$
the Lazarsfeld bundle $F_{C,A}$, which by definition is the kernel
of the evaluation map $H^0(C,A)\otimes \ko_X\to A$. Here $A$ is
viewed as a sheaf on $X$  supported on $C$. Thus, there is a short
exact sequence
$$\xymatrix{0\ar[r]&F^{\phantom{C^0}}_{C,A}\ar[r]&H^0(C,A)\otimes\ko_X\ar[r]&A\ar[r]&0}$$
and it is not difficult to see that $F_{C,A}$ really is  locally
free. Dualizing yields an exact sequence
\begin{equation}\label{eqn:sesLaz}\xymatrix{0\ar[r]&H^0(C,A)^*\otimes
\ko_X\ar[r]&F_{C,A}^{*\phantom{C^0}}\ar[r]&A^*\otimes\omega_C\ar[r]&0.}\end{equation}

The crucial result for our discussion is the following
observation.

\begin{lem}\label{lem:Laz} \cite[Lemma 1.3]{Laz}
If $|C|$ is indecomposable, i.e.\ $|C|$ does not contain any
reducible curves, then the bundle $F_{C,A}$ is simple.\qqed
\end{lem}

Clearly, the assumption on $|C|$ is satisfied if $\Pic(X)$ is
generated by $\ko(C)$ and we shall restrict to this case. So let
from now on $X$ be a complex projective K3 surface with
$\rho(X)=1$, let $H\in\Pic(X)$ be the ample generator and write
$(H.H)=2g-2$. Then the generic curve $C\in|H|$ is smooth of genus
$g$. (Indeed, by Bertini it suffices to show that $|H|$ has no
base points and according to \cite[Cor.\ 3.2]{SDonat} there are no
base points outside the fixed components which do not exist,
because $|H|$ is indecomposable.)

Now choose $r_0$ and $d$ such that
\begin{equation}\label{eqn:Lazcond} d<g+r_0\phantom{PP} {\rm
and}\phantom{PP}\rho(r_0,d,g)=0.
\end{equation} We will only need the following immediate consequence
of \cite{Laz}:

\begin{prop}\label{prop:BNgeneral}
For generic $C\in|H|$ there exists a line bundle $A\in
W^{r_0}_d(C)$ satisfying (\ref{eqn:Lazcond1}).
 \end{prop}

\begin{proof} In fact we will show that for a generic curve $C\in|H|$
any $A\in W^{r_0}_d(C)$ satisfies (\ref{eqn:Lazcond1}). Since
$\rho(r_0,d,g)=0$ and hence $W^{r_0}_d(C)\ne\varnothing$ (for any
smooth $C$), this proves the assertion.

Let $C\in |H|$ be generic and let $A\in W_d^{r_0}(C)$. We first
check $h^0(C,A)=r_0+1$. If not, then
$W_d^{r_0+1}(C)\ne\varnothing$. On the other hand, by our
assumption (\ref{eqn:Lazcond}) we have
$\rho(r_0+1,d,g)=\rho(r_0,d,g)-(g-d+r_0+1)-(r_0+2)
=0+d-g-2r_0-3<0$ and thus $W_d^{r_0+1}(C)=\varnothing$, as the
generic smooth curve in $|H|$ is Brill--Noether general according
to \cite{Laz}.

Next, for generic $C\in|H|$ any $A\in W^{r_0}_d(C)$ is globally
generated. Otherwise $W^{r_0}_{d-1}(C)\ne\varnothing$. This would
again contradict that $C$ is Brill--Noether general, for
$\rho(r_0,d-1,g)=\rho(r_0,d,g)-(r_0+1)<0$.

The calculation for $A^*\otimes\omega_C$ is similar. First observe
$A^*\otimes\omega_C\in W^{r_0-d-1+g}_{2g-2-d}(C)$ and
$h^0(C,A^*\otimes\omega_C)=r_0-d+g$ by Serre duality and
Riemann--Roch. Then, using the assumption $d<g+r_0$ in
(\ref{eqn:Lazcond}) (so far $d<g+2r_0$ was enough), one checks
$\rho(r_0-d-1+g,2g-2-d-1,g)<0$ and hence
$W^{r_0-d-1+g}_{2g-2-d-1}(C)=\varnothing$. The latter shows in
particular that $A^*\otimes\omega_C$ is globally generated.
\end{proof}

We continue to assume  $\Pic(X)=\ZZ H$. Consider a spherical
bundle $E$ on $X$ and let $v(E)=(r,kH,s)$ with $k\equiv\pm1(r)$.
By tensoring with powers of $H$ and dualizing we can modify $E$
such that $k=1$. As these operations do not affect whether
$v^\Ch(E)\in R(X)$, we will assume henceforth that $k=1$. Since
$E$ is spherical, one has $(H.H)-2rs=-2$ or, in other words,
$rs=g$.

Next we would like to relate $E$ to a particular Lazarsfeld
bundle, but a priori it is not clear that $E$ fits in a short
exact sequence of the form $0\to\ko_X^r\to E\to M\to0$  with $M$ a
line bundle on a generic $C\in|H|$ (cf.\ (\ref{eqn:sesLaz})).
However, we will see that this is possible for the right choice of
$r_0$ and $d$. To be more precise, let
$$d:=g-1-s+r\phantom{PP}{\rm and}\phantom{PP}r_0:=r-1.$$
(If the wished for exact sequence $0\to \ko_X^r\to E\to M\to0$ is
of the form (\ref{eqn:sesLaz}), then the Riemann--Roch formula
$\chi(M)=-d+g-1$ together with $\chi(M)=\chi(E)-2r=s-r$ dictates
this choice.)

A straightforward computation reveals that with this choice
$\rho(r_0,d,g)=g-rs=0$ and $d<g+r_0$. The latter is equivalent to
$s>0$ which follows from $g=rs$ and $r>0$. Thus Proposition
\ref{prop:BNgeneral} applies and we find for generic $C\in|H|$ a
line bundle $A\in W^{r_0}_d(C)$ satisfying (\ref{eqn:Lazcond}).
This then yields a short exact sequence of the form
(\ref{eqn:sesLaz}). Moreover, $F_{C,A}^*$ is simple by Lemma
\ref{lem:Laz} and $v(F_{C,A}^*)=(r,H,s)=v(E)$.

By \cite[Prop.\ 3.14]{Muk}, any spherical bundle on a K3 surface
with Picard number one is $\mu$-stable. Mukai also proves that
rigid $\mu$-stable vector bundles with given Mukai vector are
unique (see also \cite[Thm.\ 6.16]{HL}). Hence $E\cong F_{C,A}^*$.

Thus as a consequence of  \cite{Laz} we proved

\begin{cor}\label{cor:genericcurveses}
Let $E$ be a spherical bundle on a K3 surface $X$ with
$\Pic(X)=\ZZ H$ and such that $v(E)=(r,H,s)$. Then for any generic
smooth curve $C\in|H|$ there exists  a line bundle $M$ on $C$ and
a short exact sequence
$$0\to\ko_X^r\to E\to M\to0.$$\qqed
\end{cor}

\begin{remark} Corollary \ref{cor:sphsame}, which also works for
$\rho(X)>1$, shows that in any case $v^\Ch(E)=v^\Ch(F_{C,A}^*)$.
So we do not actually need $E\cong F_{C,A}^*$, but only that the
simple $F_{C,A}^*$ exists.
\end{remark}

The rough idea of the next step is to let degenerate the smooth
generic curve $C$ to a rational curve $C_0\in|H|$, which always
exist due to Mumford (cf.\ \cite{MM} or \cite{BHPV}). At the same
time, $M$ will deform to a sheaf $M_0$ supported on $C_0$. Since
the right hand side in $v^\Ch(M)=v^\Ch(E)-v^\Ch(\ko_X^r)$ stays
constant in the process, one also has
$v^\Ch(M_0)=v^\Ch(E)-v^\Ch(\ko_X^r)$. But now $M_0$ is supported
on the rational curve $C_0\subset X$ and by \cite{BV} this implies
$v^\Ch(M_0)\in R(X)$. Hence $v^\Ch(E)\in R(X)$.

This can be made rigorous as follows: Consider $Z:={\rm
Gr}(r,H^0(X,E))$ and the non-empty Zariski open subset $U\subset
Z$ of all subspaces $V\subset H^0(X,E)$ such that
$V\otimes\ko_X\to E$ is injective with cokernel $M$ being a line
bundle on a smooth curve $C$. Then $C\in|H|$, for ${\rm
c}_1(E)=H$. By Corollary \ref{cor:genericcurveses} the set $U$ is
not empty and in fact the composition
$$U\to\Pic^d(\kc/|H|)\to|H|$$ is dominant.

Here $\kc\to|H|$ is the linear system together with its universal
curve and $\Pic^d(\kc/|H|)\to|H|$ denotes the compactified
relative Jacobian variety (or Simpson's moduli space of stable
pure sheaves).

The morphism $U\to\Pic^d(\kc/|H|)$ can be compactified to a
morphism $\varphi:Z'\to\Pic^d(\kc/|H|)$ where $Z'$ is some
projective variety containing $U$ as a dense open subset. For the
following, we can assume that  the universal sheaf $\km$ on
$\Pic^d(\kc/|H|)\times X$ exists, otherwise pass to some
projective variety dominating $\Pic^d(\kc/|H|)$. Then the
pull-back $\kn:=(\varphi\times{\rm id})^*\km$ on $Z'\times X$ has
the property that $v^\Ch(\kn_t)=v^\Ch(E)-v^\Ch(\ko_X^r)$ for any
closed point $t\in U$. Hence also
$v^\Ch(\kn_{t_0})=v^\Ch(E)-v^\Ch(\ko_X^r)$ for any closed point
$t_0\in Z'$ in the boundary.

\begin{remark}\label{rem:spec}
This last argument makes use of the specialization map for Chow
groups. Consider first a family  $\kx\to S$ over a smooth
irreducible curve. Let $t\in S$ be a closed point and $\eta\in S$
be the generic point. Denote by $\kx_t$ and $\kx_\eta$ the
corresponding fibres, which we regard as varieties over $k(t)$
resp.\ $k(\eta)=K(S)$. The closure of any cycle on $\kx_\eta$
yields a cycle on $\kx$ which can then be restricted to the closed
fibre $\kx_t$. Rational equivalence is preserved in the process,
so that we get the specialization map
\begin{equation}\label{eqn:spec}
\Ch^*(\kx_\eta)\to\Ch^*(\kx_t). \end{equation} See \cite[Ch.\
20]{Fulton} for details when $S$ is the spectrum of a discrete
valuation ring with the two points $t$ and $\eta$. For an
arbitrary (smooth and irreducible) base $S$ one constructs by
recurring blow-ups (see e.g.\ \cite[II, Exer.\ 4.12]{Ha}) a
morphism $\Spec(R)\to S$, with $R$ a discrete valuation ring,
mapping the closed (resp.\ generic) point to $t$ (resp.\ $\eta$).
Then pull-back the family $\kx$ to $\Spec(R)$ and apply the
construction for discrete valuation rings. Note that by
construction for any cycle $\alpha\in\Ch^*(\kx)$ the restriction
$\alpha_t\in\Ch^*(\kx_t)$ equals the image under the
specialization map of the restriction $\alpha_\eta$, cf.\
\cite[20.3.1]{Fulton}.

This specialization technique  applies to our case, as the short
exact sequences associated to any $V\subset H^0(X,E)$ in $U$ glue
to a short exact sequence over $U$ and hence over the generic
point $\eta\in U\subset Z'$. Thus one has
$v^\Ch(\kn_\eta)=v^\Ch(E_\eta)-v^\Ch(\ko_{X_\eta}^r)$, where
$E_\eta$ is obtained by base change $X_\eta:=X\times_\CC\eta\to
X$.
\end{remark}

Since $U\to|H|$ is dominant, there exists a closed point $t_0\in
Z'$ which under $$Z'\to\Pic^d(\kc/|H|)\to|H|$$ maps to a closed
point corresponding to an irreducible rational curve $C_0\in|H|$.
Thus $\kn_{t_0}$ is supported on an irreducible rational curve.
But then $v^\Ch(\kn_{t_0})\in R(X)$ and hence $v^\Ch(E)\in R(X)$.
Thus we have proved

\begin{prop}\label{prop:mainPicone}
Let $X$ be a K3 surface with $\Pic(X)=\ZZ H$. If $E$ is a
spherical bundle with $v(E)=(r,kH,s)$ and $k\equiv\pm1(r)$, then
$v^\Ch(E)\in R(X)$.\qqed
\end{prop}

\begin{remark}
Note that the above arguments also work for $\rho(X)\geq2$
whenever $U\to |H|$ is dominant, but without using Lazarsfeld's
result this seems difficult.
\end{remark}

\begin{remark}\label{rem:Mumford}
There is an alternate argument going back to Mumford that would
replace the degeneration argument. One can show that the set of
effective cycles $Z\in{\rm S}^n(X)$ rationally equivalent to a
given one $Z_0\in{\rm S}^n(X)$ forms a countable union of
irreducible Zariski closed subsets. The countability stems from
the fact that the number of irreducible components of the Hilbert
scheme of all subvarieties is countable. See \cite[Ch.\
22]{Voisin} for an account.

By taking sections of the line bundles $\kn_t$ one obtains
effective cycles on $X$, which for all $t$ in the open subset $U$
are rationally equivalent to each other (and to ${\rm c}_2(E)$).
But then this holds for any cycle in the closure of the image of
$U\to {\rm S}^n(X)$, which necessarily contains a cycle that is
contained in a rational curve.
\end{remark}

In order to fully prove Proposition \ref{prop:sphvb}, it remains
to treat the case $\rho(X)\geq2$. There are essentially two
arguments involved:

\medskip

i) Show that on a K3 surface $X$ with $\rho(X)\geq2$ tensoring
with line bundles and dualizing brings the Mukai vector of any
spherical bundle $E$ into the form $v(E)=(r,H,s)$, where $H$ is a
primitive ample line bundle.

ii) View any polarized K3 surface $(X,H)$ with $\rho(X)\geq2$ as a
degeneration of a polarized K3 surface of Picard number one. Then
use Proposition \ref{prop:mainPicone} and a degeneration argument.

\medskip

The first question is purely numerical: Suppose $E$ is a spherical
vector bundle with $v(E)=(r,{\rm c}_1(L),s)$. Write ${\rm
c}_1(L)=k\ell$ for some primitive $\ell\in{\rm NS}(X)$ and
$k\in\ZZ$. As $E$ is spherical, one has $k^2(\ell.\ell)-2rs=-2$.
Thus $k$ and $r$ are coprime, for $(\ell.\ell)$ is even. Assuming
$\rho(X)\geq2$, there exists a line bundle $M\in\Pic(X)$ such that
${\rm c}_1(E\otimes M)=k\ell+r{\rm c}_1(M)$ is primitive and
ample. (Indeed, complete $e_1:=\ell$ to a basis $e_1,e_2,\ldots,
e_\rho$ of ${\rm NS}(X)$ and choose $M$ such that ${\rm
c}_1(M)=\sum a_ie_i$ with $a_2=\pm r^n$. Then $k+ra_1$ and
$ra_2=\pm r^{n+1}$ are coprime and for $n\gg0$ the coefficients
$a_1,a_3,\ldots, a_\rho$ can be chosen such that $ke_1+r{\rm
c}_1(M)$ is contained in the ample cone, which is open.)

Since $v^\Ch(E)\in R(X)$ is equivalent to $v^\Ch(E\otimes M)\in
R(X)$, it suffices to consider the following situation: $(X,H)$ is
a polarized K3 surface with $H$ primitive and $E$ is a spherical
vector bundle on $X$ with $\det(E)=H$.

In step ii) we choose a smooth projective family of polarized K3
surfaces $\pi:(\kx,\kh)\to D$ over a curve $D$, such that a
distinguished fibre, say over the closed point $0\in D$, is
$(X,H)$, i.e.\ $\kx_0\cong X$ and $\kh_0:=\kh|_{\kx_0}\cong H$,
and such that the general fibre has Picard number one. More
precisely, for all except countably many closed points $t\in D$
one has $\rho(\kx_t)=\ZZ\kh_t$.

The obstructions to deform the spherical bundle $E$ on the central
fibre $X$ to a bundle on the nearby fibres in $\kx\to D$ are
contained in $\Ext_X^2(E,E)$ and their traces are the obstructions
to deform $\det(E)=H$ sideways. Since $\kh$ exists, the latter
must be trivial. As $E$ is spherical, the trace free part of
$\Ext^2_X(E,E)$ is trivial and hence $E$ deforms to a vector
bundle $\ke$ on $\kx$, possibly after shrinking $D$. So
$\ke_0:=\ke|_{\kx_0}\cong E$ and, by semi-continuity (shrink $D$
again if necessary), the restriction $\ke_t$ to any other fibre
$\kx_t$ is as well spherical. A degeneration argument then yields
the following result, which completes the proof of Proposition
\ref{prop:sphvb}.

\begin{prop} Suppose $E$ is a spherical vector bundle on a K3
surface with $\rho(X)\geq2$. Then $v^\Ch(E)\in R(X)$.
\end{prop}

\begin{proof}
Again, there are two ways of proving this (cf.\ Remark
\ref{rem:spec} and the discussion following it). One can argue as
Mumford and say that either $v^\Ch(\ke_t)\in R(\kx_t)$ for all
closed points $t\in D$ or for only a countable number of them
(cf.\ Remark \ref{rem:Mumford} which one easily adapts to the
relative setting). Since over $\CC$ the number of closed points
$t\in D$ with $\rho(\kx_t)=1$ is uncountable and for them
$v^\Ch(\ke_t)\in R(\kx_t)$ by Proposition \ref{prop:mainPicone},
we must have $v^\Ch(\ke_t)\in R(\kx_t)$ for all closed points
$t\in D$ and in particular for $t=0$. Hence, $v^\Ch(E)\in R(X)$.

An alternate argument would be the following. Consider the
relative Grassmannian ${\rm Gr}(r,\pi_*\ke)$ with fibres ${\rm
Gr}(r, H^0(\kx_t,\ke_t))$ (at least over a non-empty open subset
of $D$ to which we tacitly restrict). Then let $\ku\subset {\rm
Gr}(r,\pi_*\ke)$ be the open subset of subspaces $V\subset
H^0(\kx_t,\ke_t)$ inducing short exact sequences of the form $0\to
V\otimes\ko_{\kx_t}\to\ke_t\to M\to0$ with $M$ a line bundle on
some smooth curve on $\kx_t$ in the ample linear system $|\kh_t|$.
As explained earlier, if $\rho(\kx_t)=1$, the natural morphism
$\varphi_t:\ku_t\to|\kh_t|$ is dominant. But this is an open
condition. Hence $\varphi_t$ is actually surjective on a Zariski
open subset of closed points $t\in D$ and thus over the generic
point $\eta\in D$. Then imitate the degeneration argument for the
linear system on $\kx_\eta$, that shows that $v^\Ch(\ke_\eta)\in
R(\kx_\eta)$. Now use the specialization map
$\Ch(\kx_\eta)\to\Ch(\kx_0)$.
\end{proof}

\section{K3 surfaces over number fields}\label{sect:K3overnumber}

The situation changes dramatically if instead of K3 surfaces over
$\CC$ one considers smooth projective K3 surfaces defined over a
number field or over $\bar\QQ$. In fact, a general conjecture of
Beilinson and Bloch (see \cite{BlochCrelle,RSS}) applied to this
case can be stated as follows:

\begin{conj} \label{conj:BB}If $X$ is a smooth projective K3 surface
over a number field $K$ or $\bar\QQ$, then $${\rm
deg}:\Ch^2(X)\otimes\QQ\congpf\QQ.$$
\end{conj}

How does this compare to \cite{BV} and to the results of the
previous sections? Choose an embedding $K\subset \CC$ and let
$X_\CC:=X\times_K\CC$ be the induced complex K3 surface. A
folklore argument shows that for arbitrary $X$ the kernel of the
natural map
$$\Ch^*(X)\to\Ch^*(X_\CC)$$ is torsion.
 Thus,
Conjecture \ref{conj:BB} can be rephrased as

\begin{conj}{\rm\bf (Bloch--Beilinson for K3 surfaces)}\label{conj:BB2} If $X$ is a smooth projective K3
surface over a number field $K\subset\CC$ (or $\bar\QQ$), then the
pull-back yields an injection
$$\Ch^*(X)\otimes\QQ~\hookrightarrow~ R(X_\CC)\otimes\QQ\subset\Ch^*(X_\CC)_\QQ.$$
\end{conj}

To prove the conjecture, it suffices to show that any $K$-rational
point $x\in X(K)$ satisfies $[x]=c_X\in\Ch^2(X_\CC)$. This would
follow from the a priori stronger statement that any $K$-rational
point $x\in X(K)$ lies on a rational curve which is called a
`logical possibility'  by Bogomolov, cf.\ \cite{BT}. Note that
there are other classes of points on K3 surfaces which are known
to have fundamental class in $R(X_\CC)$, e.g.\ in
 \cite{ML} this is shown for points that can be written as sums of
 a torsion point on an elliptic curve and a point in the
 intersection of the elliptic curve with a rational curve.

From the derived point of view, any $K$-rational point $x$ defines
a \emph{semi-rigid object} $$k(x)\in\Db(X),$$ where $\Db(X)$ is
viewed as a $K$-linear triangulated category. By definition, an
object $E\in\Db(X)$ is called semi-rigid if $\Ext^*_X(E,E)\cong
H^*({\rm S}^1\times {\rm S}^1,K)$. For comparison, recall that $E$
was called spherical if $\Ext^*_X(E,E)\cong H^*({\rm S}^2,K)$.

The techniques of this article do not allow to treat semi-rigid
objects, but they do show that their simpler spherical cousins
behave as expected.

\begin{prop}
Let $E\in\Db(X)$ be a spherical object on a smooth projective K3
surface $X$ over a number field $K\subset\CC$ (or  $\bar\QQ$) such
that $\rho(X_\CC)\geq2$. Then under $\Ch^*(X)\to\Ch^*(X_\CC)$ its
Mukai vector $v^\Ch(E)\in\Ch^*(X)$ is mapped to $R(X_\CC)$.
\end{prop}

\begin{proof}
This is an immediate consequence of Corollary  \ref{cor:sphobR}.
Indeed, flat base change turns $E$ into a spherical object
$E_\CC\in\Db(X_\CC)$ whose Mukai vector is contained in
$R(X_\CC)$.
\end{proof}


By means of the proposition one can now produce non-trivial
classes on K3 surfaces over number fields that are  contained in
$R(X_\CC)$. In other words, these classes behave as predicted by
Conjecture \ref{conj:BB2}, but for a less geometric reason than
e.g.\ rational points contained in rational curves.

 As it turns out, in fact all spherical objects on the
complex K3 surface $X_\CC$, which although rigid exist in
abundance, are defined over $\bar\QQ$. This is

\begin{prop}
 Let $X$ be a smooth projective K3 surface over a number
field $K\subset \CC$. Then any spherical object $F\in\Db(X_\CC)$
is defined over some finite extension $L/K$, i.e.\ there exists a
spherical object $E\in\Db(X_L)$ such that $E_\CC\cong F$.
\end{prop}

\begin{proof}
This uses a standard argument that roughly says that all points of
a zero-dimensional moduli space representing a moduli functor
defined over an algebraically closed field are defined over the
same field. E.g.\ any line bundle on $X_\CC$ is defined over
$\bar\QQ$ (and hence over some finite extension of $K$), because
the Picard variety for $X_{\bar\QQ}$ lives over $\bar\QQ$.

In our case we use Inaba's moduli space of simple complexes (cf.\
\cite{Lieb} for a more general setting). Consider the functor
${\rm\bf Splcx}_{X_{\bar\QQ}}$ on the category of locally
noetherian schemes over $\bar\QQ$, which in particular sends the
spectrum of any finitely generated field extension $L/\bar \QQ$ to
the set of isomorphism classes of all bounded complexes
$E\in\Db(X_L)$ with $\Ext^0_{X_L}(E,E)=L$ and
$\Ext^i_{X_L}(E,E)=0$ for $i<0$. Then it is shown in \cite{Inaba1}
that the \'etale sheafification ${\rm \bf
Splcx}^{et}_{X_{\bar\QQ}}$ is represented by an algebraic space
over $\bar\QQ$, which we denote ${\rm\bf Spl}$.

Any bounded complex $F\in\Db(X_\CC)$ is defined over some finitely
generated field extension $L/\bar\QQ$. Thus, a spherical
$F\in\Db(X_\CC)$ can be seen as an $L$-rational point of ${\rm\bf
Spl}$. The vanishing of $\Ext^1(F,F)$ shows that the Zariski
tangent space at the corresponding point in ${\rm\bf Spl}_L$ is
trivial. In particular, locally around the point corresponding to
$F$  the algebraic space ${\rm \bf Spl}$ is zero-dimensional and
we may therefore assume it is a scheme over $\bar\QQ$.

To conclude use the following straightforward argument from
commutative algebra. Let $A$ be a finitely generated $k$-algebra
over an algebraically closed field $k$, let $L/k$ be  any
extension, and let $B:=A\otimes_kL$. Let ${\bf n}\subset B$  be a
maximal ideal and suppose that ${\bf m}:=A\cap{\bf n}$ is maximal
in $A$. If now ${\bf n}/{\bf n}^2=0$, then $k=A$. Indeed,
Nakayama's lemma immediately shows that $B$ must be a field and
hence ${\bf m}=A\cap{\bf n}=0$, i.e.\ $A$ is a field.
 Since
$A$ is a finitely generated algebra over the algebraically closed
field $k$, this yields $A=k$. In order to reduce to the case that
${\bf m}$ is maximal, i.e.\ that the $L$-rational point of
${\rm\bf Spl}$ is a closed point, take a generic closed point $P$
in the Zariski closure of the image of $\Spec(L)\to{\rm\bf Spl}$.
By semi-continuity it will correspond to a spherical object on
$X_{k(P)}$. Then the above argument applies and shows that $P$ is
isolated and therefore equals the original $L$-rational point.

This shows that any spherical object $F\in\Db(X_\CC)$ is defined
eventually over $\bar\QQ$ and hence over some finite extension of
$K$.
\end{proof}

The deformation techniques used to prove Proposition
\ref{prop:weak} would allow one to avoid moduli spaces of simple
complexes and to work solely with moduli spaces of bundles, but
the above proof seems more conceptual.

\begin{remark}
As the reader will have noticed, the proof also shows that the
Fourier--Mukai kernel $\kf$ of any autoequivalence
$\Phi_\kf:\Db(X_\CC)\congpf\Db(X_\CC)$ is defined over $\bar\QQ$.
In other words, ${\rm Aut}(\Db(X_{\bar\QQ}))\cong{\rm
Aut}(\Db(X_\CC))$. Of course, the same holds for the set of
equivalences between two different K3 surfaces both defined over
$\bar\QQ$.
\end{remark}


\end{document}